\documentclass{ifacconf}
\usepackage{color}
\usepackage{graphics}
\usepackage{graphicx} 
\usepackage{natbib}     
\usepackage{subfigure}
\usepackage{amsmath}
\usepackage{amssymb}
\usepackage{algorithm}
\usepackage{algorithmic}
\usepackage{enumitem}
\usepackage{setspace}

\setitemize[1]{itemsep=5pt,partopsep=2pt,parsep=\parskip,topsep=5pt}

\newtheorem{proposition}{Proposition}

\newcommand{\proof}{\textbf{Proof. }}
\newcommand{\Qed}{\hfill$\blacksquare$}

\begin{document}
\begin{frontmatter}
\title{Distributed Multi-Building Coordination for Demand Response\thanksref{footnoteinfo}} 
\thanks[footnoteinfo]{The first two authors contributed equally.}
\author[First]{Junyan Su} 
\author[First]{Yuning Jiang} 
\author[Second]{Altu\u{g} Bitlislio\u{g}lu}
\author[Second]{Colin N. Jones}
\author[First]{Boris Houska} 
\address[First]{School of Information Science and Technology, ShanghaiTech University, Shanghai, China\\ (email: sujy, jiangyn, borish@shanghaitech.edu.cn).}
\address[Second]{Automatic Control Laboratory, Ecole Polytechnique F\'ed\' erale de Lausanne (EPFL), Lausanne, Switzerland\\ (e-mail: altug.bitlislioglu@alumni.epfl.ch, colin.jones@epfl.ch).}

\begin{abstract}                
This paper presents a distributed optimization algorithm tailored for solving optimal control problems arising in multi-building coordination. The buildings coordinated by a grid operator, join a demand response program to balance the voltage surge by using an energy cost defined criterion. 
In order to model the hierarchical structure of the building network, we formulate a distributed convex optimization problem with separable objectives and coupled affine equality constraints. 
A variant of the Augmented Lagrangian based Alternating Direction Inexact Newton (ALADIN) method for solving the considered class of problems is then presented along with a convergence guarantee. 
To illustrate the effectiveness of the proposed method, we compare it to the Alternating Direction Method of Multipliers (ADMM) by running both an ALADIN and an ADMM based model predictive controller on a benchmark case study.
\end{abstract}

\end{frontmatter}

\section{Introduction}
Energy generation is undergoing a rapid transition from fossil fuels to renewable sources~(\cite{Liserre2010}), which poses a challenge to balance the unpredictable generation demand due to the highly stochastic nature of renewable energy sources, and requiring advanced ancillary service providers. Recently, Demand Response (DR) programs utilizing the flexibility of power demand to provide services have been considered in the power systems community~(\cite{Siano2014}). These programs cover collective load shifting, real time power regulation for load balancing and capacity firming, which has been applied to mitigate the uncertainty in renewable power generation effectively~(\cite{Bitlislioglu2018}). Because commercial buildings, which are equipped with available heating, ventilation and air conditioning (HVAC) systems, have a potential to collectively offer ancillary services to the power grid~(\cite{Oldewurtel2012}). Smart grids connecting multiple commercial buildings were developed recently in the DR program to match the increasing power scale. In this setting, individual buildings are coupled via the grid operator. This yields a coordination problem, which can be put in the generic framework of multi-agent optimization and control~(\cite{Bitlislioglu2018}).  

Typically, in order to meet the real-time requirement, these multi-agent coordination problems are embedded in a Model Predictive Control (MPC) framework~(\cite{Rawlings2017}), where the resulting problems have to be solved once during each sampling time, which requires an efficient online solver. 
For this purpose, distributed algorithms have already been developed~(\cite{Bitlisliouglu2017b,Boyd2011,Braun2018}). A class of these approaches is based on decomposition methods, including primal and dual decomposition. In~(\cite{Rantzer2009,Richter2011}), gradient-based dual decomposition methods are used to solve the concave dual problem. Alternatively, semi-smooth Newton methods~(\cite{Frasch2015}) can be applied.
In~(\cite{Bitlisliouglu2017b}), an interior point method based on primal decomposition was proposed, which writes all the inequality constraints into the objective by using a primal barrier function~(\cite{Boyd2004}). As a follow-up,~(\cite{Bitlisliouglu2017c}) proposed a primal-dual interior point method, which further decomposes the resulting Newton-step. However, such Newton-type methods are in general only convergent if they are equipped with additional smoothing heuristics and line-search routines.
Compared to the decomposition method, the Alternating Direction Method of Multipliers (ADMM) has more reliable convergence properties~(\cite{Boyd2011,Hong2017}). Many variants of ADMM exploit the hierarchical structure~(\cite{Boyd2011,Goldstein2014}), but, in practice, a heuristic pre-conditioner is required to enhance convergence~(\cite{Donoghue2016}).

This paper considers the case that the grid operator coordinates a group of commercial buildings, which joins a DR program. Section~2 introduces the problem formulation based on~(\cite{Bitlisliouglu2017b}). Then, we reformulate the original problem by exploiting the decomposed structure of the building network, in which the local variables are hidden. A distributed optimization problem is thus yielded, where the decoupled objectives are non-smooth Piece-Wise Quadratic (PWQ) functions with linear coupling constraints. For solving this problem in the context of MPC, Section~3 proposes a tailored Augmented Lagrangian based Alternating Direction Inexact Newton (ALADIN) method~(\cite{Houska2016}), which comes along with a convergence guarantee. ALADIN recently has been proposed to solve multi-agent optimization problems~(\cite{Jiang2017a,Engelmann2019}). Similar to ADMM,
it requires the local agents to solve small-scale decoupled problems and the central entity to deal with a linear equation in every iteration. For this variant, a warm-start strategy is further proposed to improve its online performance. As a result, compared to ADMM, ALADIN takes the same communication effort per iteration while requiring much fewer iterations to converge to a desired accuracy. This is illustrated by running both an ALADIN and an ADMM based MPC controller in a benchmark case study.

\textit{Notation:} The set of symmetric, positive (semi-)definite matrices in $\mathbb R^{n \times n}$ is denoted by $(\mathbb S_{+}^{n})\,\mathbb S_{++}^{n}$. We use notation $\mathbf{1}_n = [1\,\ldots\,1]^\top\in\mathbb{R}^{n}$ for all $n\in\mathbb{N}$. For a given matrix $\Sigma \in \mathbb S_{+}^{n}$ the notation
\[
\left\| x \right\|_{\Sigma} = \sqrt{ x^\top \Sigma x }
\]
is used. Moreover, we call a function $f: \mathbb R^{n} \to \mathbb R \cup \{ \infty \}$ strongly convex with $\Sigma \in \mathbb S_{+}^n$, if the inequality
\[
f(t x + (1-t)y) \leq t f(x) + (1-t)f(y) - \frac{1}{2} t (1-t) \left\| x-y \right\|_{\Sigma}^2
\]
is satisfied for all $x,y \in \mathbb R^{n}$ and all $t \in [0,1]$. Finally, the Kronecker product of two matrices $A \in \mathbb{R}^{k \times l}$ and $B \in \mathbb{R}^{m \times n}$ is given by $$A \otimes B = \left(a_{ij}B\right)_{i,j} \in \mathbb{R}^{km \times ln}.$$

\section{Problem Formulation}
This section introduces a hierarchical optimal control problem for coordinating a group of commercial buildings, which join a Demand Response (DR) program. 

\subsection{Tracking Model of Single Building}
This paper considers a building network in which each building has a central heating control system. The dynamics of the $i$-th building can be described by a linear input-output system, 
\begin{equation}\label{eq::dyn}
\begin{split}
x_{i,k+1} \;=\;& A_ix_{i,k} + B_iu_{i,k} + w_{i,k}\;,\\
y_{i,k} \; =\; & C_ix_{i,k} + D_iu_{i,k}\;,
\end{split}
\end{equation}
where state $x_{i,k}\in\mathbb{R}^{n_{x_i}}$ denotes the temperatures of the thermal zones of the $i$-th building, $u_{i,k}\in\mathbb{R}^{n_{u_i}}$ the thermal cooling energy input to the building at time step $k$ and $w_{i,k}\in\mathbb{R}^{n_{w_{i}}}$ the system disturbance. The output $y_{i,k}\in\mathbb{R}^{n_{y_i}}$ denotes the mean zone temperatures and coefficient matrices $A_i$, $B_i$, $C_i$, $D_i$ depend on the specific buildings.

For a given reference room temperature $y_i^\mathrm{ref}$, the following tracking optimal control problem can be constructed,
\begin{equation}\label{eq::LocalProb}
\begin{split}
\min_{u_i}\quad& \sum_{k=0}^{N-1}\left(\left\| y_{i,k}-y_{i}^\mathrm{ref}\right\|_{Q_i}^2 +\left\|u_{i,k}\right\|_{R_i}^2\right)\\
\text{s.t.}\quad&\left\{
\begin{array}{rcl}
\forall\; k&\in&\{0,...,N-1\}\\[0.12cm]
x_{i,k+1} &=& A_ix_{i,k} + B_iu_{i,k} + w_{i,k},\\[0.12cm]
y_{i,k} &= & C_ix_{i,k} + D_iu_{i,k}\\[0.12cm]
x_{i,0}&=&\hat{x}_{i}, \\[0.12cm]
\underline{y}_i &\leq& y_{i,k}\; \leq\; \overline{y}_i,\\[0.12cm]
\underline{u}_i &\leq& u_{i,k} \;\leq\; \overline{u}_i
\end{array}
\right. 
\end{split}
\end{equation}
with $u_i= [u_{i,0}^\top\;\ldots\;u_{i,N-1}^\top]^\top$. Here $\hat{x}_i$ is the initial state, $[\underline{y}_i,\overline{y}_i]$ and $[\underline{u}_i,\overline{u}_i]$ denote box constraints on the system outputs and control inputs, respectively. The matrices $Q_i\in\mathbb{S}_{+}^{n_{x_i}}$ and $R_i\in\mathbb{S}_{++}^{n_{u_i}}$ are symmetric positive semi-definite and positive definite such that Problem~\eqref{eq::LocalProb} is a strongly convex quadratic programming (QP) problem~(\cite{Borrelli2003}). In the following, we represent the states and outputs as an affine function of the initial state $\hat{x}_i$ and control inputs $u_i$ by using the recursive derivation
\[
x_{i,k+1} = A_i^{k+1}\hat{x}_{i} + \sum_{j=0}^{k}A_i^{k-j}\left(B_iu_{i,j} + w_{i,k}\right) 
\]  
for $k=0,...,N-1$ such that the dynamics~\eqref{eq::dyn} can be written in the dense form
\[
\begin{split}
x_i &\;=\; \mathcal{A}_{i}\hat{x}_i + \mathcal{B}_i^u u_i + \mathcal{B}_i^ww_i\\
y_i & \;=\; \mathcal{C}_ix_i + \mathcal{D}_i u_i
\end{split}
\]
with $x_i$, $y_i$ and $w_i$ defined analogously to $u_i$. Thus, we can rewrite the objective of~\eqref{eq::LocalProb} as a quadratic cost
\[
f_i(u_i) = u_i^\top H_i u_i + 2h_i^\top u_i
\]
with matrix $H_i\succ 0$ and vector $h_i$ given by 
\[
\begin{split}
H_i &\;=\; (\mathcal{C}_i\mathcal{B}_i^u + \mathcal{D}_i)^\top\mathcal{Q}_i (\mathcal{C}_i\mathcal{B}_i^u + \mathcal{D}_i) + \mathcal{R}_i\;,\\
h_i&\;=\;H_i(\mathcal{C}_i\mathcal{A}_{i}\hat{x}_i + \mathcal{C}_i\mathcal{B}_i^ww_i - \mathbf{1}_N\otimes y_i^\mathrm{ref})\;,\\
\end{split}
\]
and the constraints as a polyhedral set 
\[
u_i\in\mathbb{U}_i:= \{u\in\mathbb{R}^{Nn_x} \mid E_i u_i \leq e_i \}
\] 
with matrix $E_i$ and vector $e_i$ given by
\[
E_i= \begin{bmatrix}
\mathcal{C}_i\mathcal{B}_i^u+\mathcal{D}_i\\
-\mathcal{C}_i\mathcal{B}_i^u-\mathcal{D}_i\\
I\\
-I
\end{bmatrix},
e_i = \begin{bmatrix}
\mathbf{1}_N\otimes  \overline{y}_i-\mathcal{C}_i\mathcal{A}_i\hat{x_i}-\mathcal{C}_i\mathcal{B}_i^ww_i\\
\mathcal{C}_i\mathcal{A}_i\hat{x_i}+\mathcal{C}_i\mathcal{B}_i^ww_i-\mathbf{1}_N\otimes\underline{y}_i\\
\mathbf{1}_N\otimes\overline{u}_i\\
-\mathbf{1}_N\otimes\underline{u}_i
\end{bmatrix}.
\] 
Here, we use notation $\mathcal{Q}_i=\text{diag}(Q_i,...,Q_i)$ and $\mathcal{R}_i = \text{diag}(R_i,...,R_i)$. The constraint set $\mathbb{U}_i$ is convex and compact~\cite{Braun2018} such that Problem~\eqref{eq::LocalProb} has a unique solution with respect to a given initial state $\hat{x}_i$. Next, we will present the coordination problem based on~\eqref{eq::LocalProb}. 

\subsection{Multi-building Coordination}
In this paper, we consider a smart grid with $M$ commercial buildings, which join a DR program and are coordinated by a grid operator. The goal of coordination is to balance the voltage surge caused by a large renewable energy generation such as solar plant. We denote by $\theta_{i,k}$ the active power injection or consumption from the $i$-th building.
It is a linear function of $u_{i,k}$, 
\begin{equation}\label{eq::theta}
\theta_{i,k} = F_{i} u_{i,k}\;,
\end{equation}
with energy transfer matrix $F_i$ for all $k=0,1,...,N-1$. Moreover, the overall power magnitudes $v_{k}\in\mathbb{R}$ can be represented as an affine map of $\theta_{i,k}$, $i=1,...,M$ given by
\begin{equation}\label{eq::coupling}
v_{k} = \sum_{i=1}^{M}G_{i}\theta_{i,k} + \tilde{v}_k
\end{equation}
for all $k=0,...,N-1$. Here, $\tilde{v}_k$ is given by the predicted voltage magnitudes at time step $k$, in p.u., $G_i$ models the effect of the injections from the $i$-th building on the overall voltage magnitudes. In addition, $v_k$ needs to satisfy the box constraints 
\begin{equation}\label{eq::boxV}
\underline{v}\;\leq\; v_k\;\leq\; \overline{v}
\end{equation}
for all $k=0,...,N-1$. The optimal coordination problem can be formulated as
\begin{equation}\label{eq::GlobalProb}
\begin{split}
\min_{p,\theta,u}\quad &\sum_{i=1}^{M}\left(f_i(u_i) + \sum_{k=0}^{N-1} \pi_k \theta_{i,k} \right)\\[0.12cm]
\text{s.t.}\quad & \left\{
\begin{array}{rcl}
\forall\;\; k&\in&\{0,1,...,N-1 \}\\[0.12cm]
\theta_{i,k} &=& F_{i} u_{i,k}\;,\;i=1,...,M\\[0.12cm]
v_{k} &=& \sum_{i=1}^{M}G_{i}\theta_{i,k} + \tilde{v}_k\;,\\[0.12cm]
\underline{v}&\leq& v_k\;\leq\; \overline{v}\;,\\[0.12cm]
u_i&\in&\mathbb{U}_i\;,\;i=1,...,M\;.
\end{array}
\right.
\end{split}
\end{equation}
The used demand response criterion is defined by 
the second term in the objective, which is the economic cost of the electricity. Here, $\pi_k$ denotes the price of the electricity. If we substitute~\eqref{eq::theta} into the objective and the constraints,~\eqref{eq::GlobalProb} becomes 
\begin{equation}\label{eq::GlobalProb2}
\begin{split}
\min_{u,v}\quad &\sum_{i=1}^{M}\left(f_i(u_i) + \sum_{k=0}^{N-1} \pi_k F_{i} u_{i,k} \right)\\[0.12cm]
\text{s.t.}\quad & \left\{
\begin{array}{rcl}
\forall\;\; k&\in&\{0,...,N-1 \}\\[0.12cm]
v_{k} &=& \sum_{i=1}^{M}G_{i}F_{i} u_{i,k} + \tilde{v}_k\;,\\[0.12cm]
\underline{v}&\leq& v_k\;\leq\; \overline{v}\\[0.12cm]
u_i&\in&\mathbb{U}_i\;,\; i=1,...,M\;.
\end{array}
\right.
\end{split}
\end{equation}
The coupling between the buildings is modeled by the global variable $v=[v_0^\top\;\cdots\;v_{N-1}^\top]^\top$ in~\eqref{eq::coupling}. In order to design an efficient distributed optimization algorithm, in the following, we will eliminate this variable and reformulate~\eqref{eq::GlobalProb2} into the standard distributed form, which is only with local variables coupled by an affine equality.

\subsection{Reformulation}

Let us introduce auxiliary variables $s_{i,k}\in\mathbb{R}^{2}$ for $i=1,...,M$, $k=0,...,N-1$. The constraints~\eqref{eq::coupling} and~\eqref{eq::boxV} can then be reformulated as inequality constraints 
\begin{equation}\label{eq::deCons}
	g_{i}(u_i,s_i)\leq 0
\end{equation}
with
\begin{align}\notag
g_i(u_{i},s_{i})=\left\{
\left(
\begin{array}{c}
\frac{\tilde{v}_k-\overline{v}}{M} +G_{i}F_{i} u_{i,k}  \\[0.16cm]
\frac{\underline{v}-\tilde{v}_k}{M} -G_{i}F_{i} u_{i,k} 
\end{array}
\right)
- s_{i,k}
\right\}_{k\in\{0,...,N-1\}}
\end{align}
for all $i=1,...,M$ and the equality affine constraints
\begin{equation}\label{eq::coupling2}
\sum_{i=1}^{M}s_i = 0 
\end{equation}
with $s_i = [s_{i,0}^\top\;\ldots\;s_{i,N-1}^\top]^\top$.
As a result, Problem~\eqref{eq::GlobalProb2} can be rewritten as 
\begin{equation}\label{eq::GlobalProb3}
\min_{z}\quad \sum_{i=1}^{M}\tilde{f}(u_i)\quad
\text{s.t.}\left\{
\begin{array}{rcl}
0&=&\sum_{i=1}^{M}s_i\quad \mid \lambda \\[0.16cm]
z_i&\in&\mathbb{Z}_i\;, \;i=1,...,M 
\end{array}
\right.
\end{equation}
with stacked variables $z_i=[u_i^\top\;s_i^\top]^\top$. Here, $\lambda$ denotes the Lagrangian multipliers of the affine equality constraints. The decoupled objectives are given by
\[
\tilde{f}_i(u_i)\;=\;f_i(u_i) + \sum_{k=0}^{N-1} \pi_k F_{i} u_{i,k} 
\]
and the constraint sets are denoted by
\[
\mathbb{Z}_i =
\left\{
\left. 
\left[
\begin{array}{c}
u_i\\[0.12cm]
s_i
\end{array}
\right]
\in\mathbb{R}^{Nn_{u_i}+2N}
\right| \begin{array}{c}
u_i\in\mathbb{U}_i\;,\;g_i(u_i,s_i)\leq 0
\end{array}
\right\}.
\]
Since~\eqref{eq::deCons} are affine inequalities, sets $\mathbb{Z}_i$ are convex polytopes. In the following, the equivalence between~\eqref{eq::GlobalProb2} and~\eqref{eq::GlobalProb3} is established. 
\begin{proposition}\label{prop::equivalence}
	If Problem~\eqref{eq::GlobalProb2} is feasible with solution $u^*$, Problem~\eqref{eq::GlobalProb3} has a solution $\hat{z}=(\hat{u},\hat{s})$ with $\hat{u}=u^*$. 
	Reversely, if Problem~\eqref{eq::GlobalProb3} is feasible with solution $\hat{z}=(\hat{u},\hat{s})$, $\hat{u}$ is the minimizer of Problem~\eqref{eq::GlobalProb2}.  
\end{proposition}
\proof
Let $u^*$ be a minimizers of Problem~\eqref{eq::GlobalProb2}. Then, we can construct a feasible point $\hat{z}=(u^*,\hat{s})$ of~\eqref{eq::GlobalProb3} as
\begin{equation}\label{eq::si}
\hat{s}_{i,k} = 
\left(
\begin{split}
& G_{i}F_{i} \hat{u}_{i,k} - \frac{1}{M}\sum_{i=1}^{M} G_i F_i \hat{u}_{i,k} \\[0.12cm]
-& G_{i}F_{i} \hat{u}_{i,k} + \frac{1}{M}\sum_{i=1}^{M} G_i F_i \hat{u}_{i,k}
\end{split}
\right)\;.
\end{equation}
For any feasible point $z = (u,s)$ of~\eqref{eq::GlobalProb3}, $u$ is also feasible for~\eqref{eq::GlobalProb2}. Thus, we have 
\[
\sum_{i=1}^M \tilde{f}(u_i^*) \leq \sum_{i=1}^M \tilde{f}({u}_i)\;.
\]
This shows that $\hat{z}$ is a minimizer of~\eqref{eq::GlobalProb3}.
Similarly, for the other direction, let $\hat{z} = (\hat{u}, \hat{s})$ be a minimizer of~\eqref{eq::GlobalProb3}. For any feasible point $u$ of~\eqref{eq::GlobalProb2}, we can construct a feasible point $z = (u,s)$ of~\eqref{eq::GlobalProb3} based on~\eqref{eq::si} such that 
\[
\sum_{i=1}^M \tilde{f}(\hat{u}_i) \leq \sum_{i=1}^M \tilde{f}(u_i)\;.
\]
Therefore, $\hat{u}$ is a minimizer of Problem~\eqref{eq::GlobalProb2}.
\Qed	

Concerning Problem~\eqref{eq::GlobalProb3}, due to the strong convexity of $f_i(\cdot)$ and compact polyhedron $\mathbb{Z}$, the optimal solution $(u^*,s^*)$ of~\eqref{eq::GlobalProb3} is unique with respect to $u^*$ but not $s^*$.
Therefore, we further introduce a least-squares regularization of $s_i$ in the decoupled objective, 
\begin{equation}\label{eq::augObj}
\mathcal{F}_i(z_i)= \tilde{f}_i(u_i) + \mu\|s_i\|_2^2 \
\end{equation}
with a sufficiently small constant $\mu>0$. This regularization enforces strong convexity of the problem and thus, uniqueness of $s_i$. Note that in practice, this small regularization does not lead to large changes of the optimal solution, which will be numerically illustrated later.  

Next, we introduce the function $\Psi_i:\mathbb{R}^{2N}\rightarrow\mathbb{R}_{\geq}$ given by the following multi-parametric QP (mpQP) problem
\begin{equation}\label{eq::Psi}
\Psi_i(s_i) \,=\, \min_{u_i}\;\;\mathcal{F}_i(z_i)\quad\text{s.t.}\;\; z_i\in\mathbb{Z}_i
\end{equation}
for all $i=1,...,M$.
According to Theorem~2 in~(\cite{Alessio2009}), $\Psi_i$ is a strongly convex PWQ function and its solution map $u_i^\star:\mathbb{R}^{2N}\rightarrow \mathbb{R}^{Nn_{u_i}}$ is Piece-Wise Affine (PWA). 
In the following, we use the notation 
\[
\Psi_i(s_i) = \frac{1}{2}s_i^\top \mathcal{S}_i(s_i) s_i + r_i(s_i)^\top s_i \;,
\]
where matrix $\mathcal{S}_i(\cdot)$ and vector $r_i(\cdot)$ are piece-wise constant with respect to $s_i$ inside different critical regions~(\cite{Borrelli2003}). For a given $s_i$, we denote the active constraints at $u_i^\star(s_i)$ by 
\[
[P_{i,1}(s_i)\;\;P_{i,2}(s_i)]\begin{bmatrix}
u_i^\star(s_i)\\ s_i
\end{bmatrix}+ p_i(s_i) = 0
\] 
with active Jacobian $[P_{i,1}(s_i)\;\;P_{i,2}(s_i)]$, the matrices $\mathcal{S}_i(s_i)$ are given by 
\begin{equation}\label{eq::Si}
\begin{split}
\mathcal{S}_i(s_i) 
=& \;\mu I+\\
&P_{i,2}(s_i)^\top [ P_{i,1}(s_i) H_i^{-1} P_{i,1}(s_i)^\top]^{-1}P_{i,2}(s_i)
\end{split}
\end{equation}
with $\mathcal{S}_i(s_i)\in\mathbb{S}_{++}^{2N}$.
Accordingly, the multi-building coordination problem can be written as
\begin{equation}\label{eq::Prob}
\min_{s}\;\; \sum_{i=1}^{M}\Psi_i(s_i)\quad
\text{s.t.}\;\;\sum_{i=1}^{M}s_i=0\;,\quad \mid \lambda \\[0.12cm]
\end{equation}
which is a strongly convex but non-smooth optimization problem. In the following, we will design an algorithm for solving~\eqref{eq::Prob} in a distributed manner.

\section{Algorithm}
This section proposes a distributed algorithm based on ALADIN~(\cite{Houska2016}) for multi-building coordination.

\begin{algorithm}[htbp!]
	\setstretch{1.08}
	\caption{ALADIN for solving~\eqref{eq::GlobalProb3}}
	\textbf{Initialization:} Initial guess $(s,\lambda)$, choose symmetric scaling matrices $\Sigma_i\succ 0$ and terminal tolerance $\epsilon>0$. \\[0.16cm]
	\textbf{Repeat:}
	\begin{itemize}
		\item[1)] Each building solves the decoupled QP
		\begin{equation}\label{eq::deProb}
			\begin{split}
				\min_{\xi_i}\;\;&\mathcal{F}_i(\xi_i) +\lambda^\top \xi_i^s+\frac{1}{2}\left\|\xi_i^s-s_i\right\|_{\Sigma_i}^2\\
				\text{s.t.}\;\;&\xi_i=(\xi_i^u,\xi_i^s)\in\mathbb{Z}_i
			\end{split}
		\end{equation}
		for $i=1,...,M$ in parallel and send solution $\xi_i^s$ to the grid operator. 
		\item[2)] Terminate if $\|s-\xi^s\| \leq \epsilon$. 
		\item[3)] The grid operator collects $\xi_i^s$ and solve the equality constrained QP
		\begin{equation}\label{eq::coQP}
			\begin{split}
				\min_{s^+}&\;\;\sum_{i=1}^{M}\frac{1}{2}\left\|s_i^+-2\xi_i^s+s_i\right\|_{\Sigma_i}^2 \\
				\text{s.t.}&\;\;\sum_{i=1}^{M}s_i^+ =0 \quad \mid \Delta \lambda^+
			\end{split}
		\end{equation}
		Then, update $\lambda^+ = \lambda + \Delta \lambda^+$ and spread $(s_i^+,\lambda^+)$ to $i$-th building for all $i=1,...,M$.
	\end{itemize}
	\label{alg::aladin}
\end{algorithm}

\subsection{Distributed Optimization with ALADIN}
Algorithm~\ref{alg::aladin} outlines to solve~\eqref{eq::Prob}
by using a tailored ALADIN algorithm. 
Similar to the standard ALADIN method, Algorithm~\ref{alg::aladin} alternates between solving~\eqref{eq::deProb} in parallel and dealing with the equality constrained QP problem~\eqref{eq::coQP} for consensus. Here, Problems~\eqref{eq::deProb} are equivalent to
\begin{equation}\label{eq::deProb2}
\min_{\xi_i^s}\;\;\Psi_i(\xi_i^s) + \lambda^\top \xi_i^s+\frac{1}{2}\left\|\xi_i^s-s_i\right\|_{\Sigma_i}^2\;,
\end{equation}
which are also mpQPs with input parameters $(\lambda,s_i)$. Thus,
the solution maps of~\eqref{eq::deProb2}, denoted by $$\xi_i^\star(\lambda,s_i)\;,\;i=1,...,M\;,$$ 
are piece-wise affine~(\cite{Alessio2009}). Due to QP~\eqref{eq::coQP} without inequality constraints,
the solution can be worked out analytically, 
\begin{subequations}\label{eq::solQP}
	\begin{align}\label{eq::dualup}
	\Delta \lambda^+ &= 2\Lambda^{-1} \left( \sum_{i=1}^{M}\xi_i^s \right)\quad \text{with }\Lambda = \sum_{i=1}^{M}\Sigma_i^{-1} \\\label{eq::priup}
	s_i^{+} &= 2\xi_i^s-s_i - \Sigma_i^{-1}\Delta \lambda^+\;,\;i=1,...,M.
	\end{align}
\end{subequations}
Here, it is clear that the grid operator only needs to collect the local solution $\xi_i^s$ and spread $\Delta \lambda^+$ to each building. Compared to ADMM, Algorithm~\ref{alg::aladin} takes exactly the same communication effort as ADMM per iteration~(\cite{Houska2016}). 

\subsection{Convergence Analysis}
As discussed in~(\cite{Jiang2019}), the iterates of Algorithm~\ref{alg::aladin} converges globally with a linear rate. Furthermore, if the scaling matrices are chosen as 
$\Sigma_i = \mathcal{S}_i(s_i)$ with $\mathcal{S}_i(s_i) = \mathcal{S}_i(s_i^*)$ the exact Hessian of $\Psi_i$ at the optimal $s_i^*$, Algorithm~\ref{alg::aladin} can further achieve a local one-step convergence under a regularity condition~(\cite{Frasch2015}). Note that this choice of $\Sigma_i$ requires prior knowledge of the optimality of~\eqref{eq::Prob}, which is in general impractical. However, in the context of Model Predictive Control (MPC), the result of the last MPC iteration can be used to choose $\Sigma_i$ online, which has a potential to satisfy $\Sigma_i = \mathcal{S}_i(s_i^*)$, and thus, local convergence can be improved. 
\subsection{Online Implementation Details}
In order to arrive at an efficient implementation, the structure of~\eqref{eq::deProb} and~\eqref{eq::coQP} can be exploited as follows.

\textbf{Online solver:} When we apply Algorithm~\ref{alg::aladin} as an online solver for MPC, we can move the primal update~\eqref{eq::priup} into the local phases such that a  simplified version of Algorithm~\ref{alg::aladin} is given by
\begin{subequations}
\label{eq::aladin}
\begin{align}\label{eq::parallel}
\text{Parallel Step}\;\;\;& 
\left\{
\begin{array}{rcl}
\lambda^+ &=& \lambda + \Delta \lambda\,,\\[0.12cm]
s_i^+ &= &2\xi_i^s - s_i - \Sigma_i^{-1}\Delta \lambda\,,\\[0.12cm]
\xi_i^+&=& \xi_i^\star(\lambda^+,s_i^+)\,,
\end{array}
\right.
\\[0.16cm]\label{eq::consensus}
\text{Consensus Step}\;&\quad \Delta \lambda^+ = 2\Lambda^{-1}\left(\sum_{i=1}^{M}{\xi_i^{s}}^+
\right)\;.
\end{align}
\end{subequations}

\textbf{Warm-start:} In an MPC scheme, the initial guess of Algorithm~\ref{alg::aladin} can be initialized by shifting the horizon,
\[
\begin{split}
s_i&\;=\;(s_{i,1}^*,\ldots,s_{i,N-1}^*,0)\;,\;i=1,...,M\;,\\
\lambda&\;=\;(\lambda_{1}^*,\ldots,\lambda_{N-1}^*,0)\;.
\end{split}
\]
This strategy has been used in the context of an ADMM based model predictive controller for smart grids~(\cite{Braun2018}). Here, $(s^*,\lambda^*)$ denotes the optimal solution of the current MPC problem. Furthermore, for Algorithm~\ref{alg::aladin}, we can set $\Sigma_i =\mathcal{S}_i(s_i^*)$, which potentially improves the local convergence as discussed in Section~3.2. Note that in practice, matrix $\mathcal{S}_i(s_i^*)$ might be ill-conditioned such that a iterative linear equation solver is required to deal with the equality constrained QP~\eqref{eq::coQP}. 

\section{Numerical Results}
This section illustrates the effectiveness of Algorithm~\ref{alg::aladin} in the MPC scheme by comparing it to the state-of-the-art method ADMM.

In our implementation, both algorithms are executed as the online solver for MPC. And the warm-start strategy discussed in Section~3.3 is applied. Furthermore, in order to obtain a fair comparison, we implemented a pre-conditioner for ADMM by performing a modified Ruiz equilibration~(\cite{Ruiz2001}) on the decoupled constraint matrices $E_i$. 

The data used to generate the benchmark is obtained by using the \texttt{EnergyPlus} toolkit~(\cite{Energyplus}) and the thermal model of buildings are generated with the \texttt{OpenBuild} toolbox~(\cite{Gorecki2015}). The length of prediction horizon is chosen by $N = 14$ with sampling time $0.5$ hour. Here, we consider three types of buildings with different scales,
\smallskip
\begin{center}
	\renewcommand{\arraystretch}{1.5}		
	\begin{tabular}{|l||c|c|c|}
		\hline
		Type & Large & Middle & Small\\ \hline
		dimension of $z_i$ & $280$    & $98$& $70$  \\
		dimension of $e_i$ & $1036$ & $308$&$196$ \\ \hline 
	\end{tabular}
\end{center}
\smallskip
where the first row represents the number of variables and the second gives the number of inequality constraints with respect to a single building. Moreover, an interval constraint of $v(k)$ is given by $v_{k}\in[0.95,1.05]$ in p.u., and the price $\pi_k=1$ is fixed for all $k=0,...,N-1$. 
\begin{figure}[htbp!]
	\centering
	\includegraphics[width=1.075\linewidth]{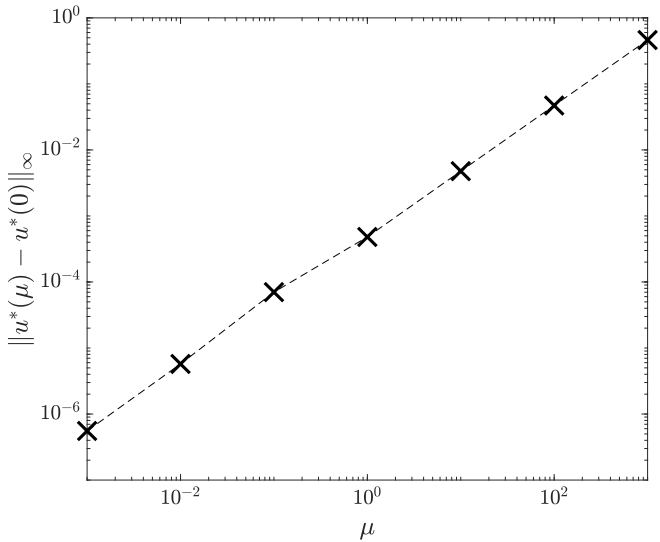}
	\caption{Difference between $u^*(\mu)$ and $u^*(0)$ in [kW].}
	\label{fig::mu}
\end{figure}

We consider a benchmark case study from~(\cite{Bitlisliouglu2017b}) with a mix of $12$ commercial buildings including $2$ Large, $7$ Middle and $3$ Small
such that there are $1456$ variables and $4816$ affine inequality constraints in total.
Regarding the choice of $\mu$, Fig.~\ref{fig::mu} illustrates the difference between the optimal solutions of~\eqref{eq::Prob} with different $\mu\neq 0$ and the result with $\mu=0$. The gap $\|u^*(\mu)-u^*(0)\|_\infty$ increases linearly with $\mu$ increasing. We set the accuracy of the online solver as $10^{-4}$ such that we choose $\mu= 10^{-1}$ in our implementation.

\begin{figure}[htbp!]
	\centering
	\includegraphics[width=1.075\linewidth]{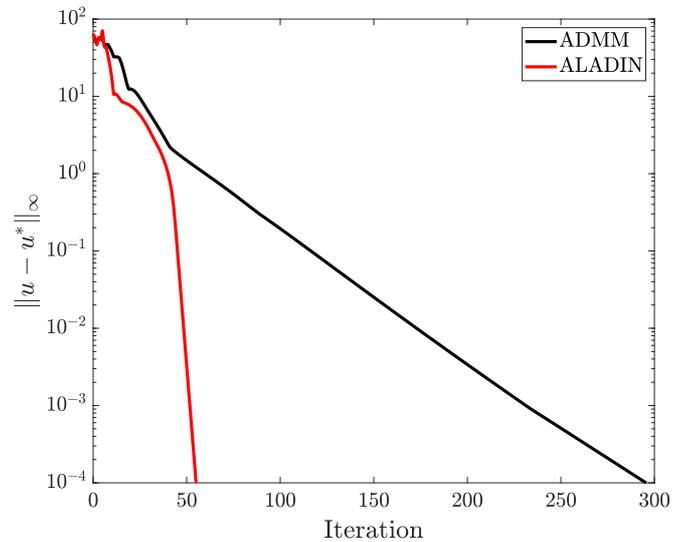}
	\caption{Convergence comparison of Case I: ALADIN vs ADMM.}
	\label{fig::convergence1}
\end{figure}

\begin{figure}[htbp!]
	\centering
	\includegraphics[width=1.075\linewidth]{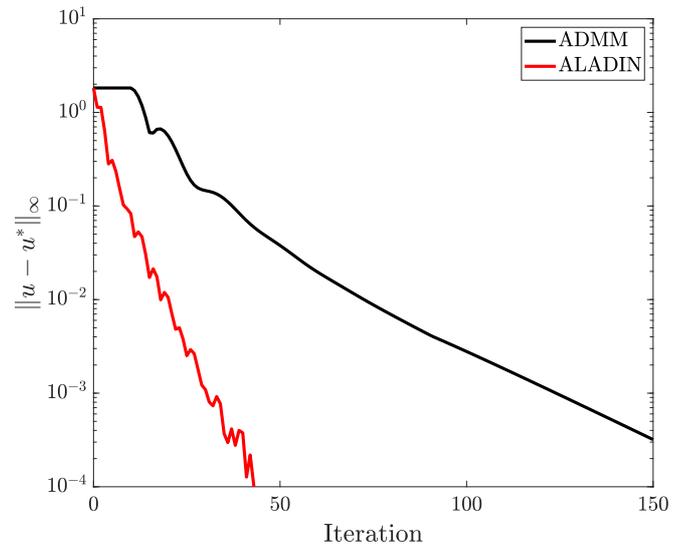}
	\caption{Convergence comparison of Case II: ALADIN vs ADMM.}
	\label{fig::convergence2}
\end{figure}

Fig.~\ref{fig::convergence1} and Fig.~\ref{fig::convergence2} show the convergence comparison between ALADIN and ADMM for two different online cases. In the first one, the optimal active set is almost the same as the previous MPC iteration. On the contrary, there are some changes of the optimal active set arising in the second case.

Fig.~\ref{fig::convergence1} shows the comparison of an example of the first case. 
For solving this particular problem, when we set the tolerance to $10^{-4}$, ALADIN is five times faster than ADMM. In this case study, $252$ constraints are active at the optimal solution. The warm-start strategy for initializing $\Sigma_i$ improves the local convergence of ALADIN. 
Fig.~\ref{fig::convergence2} shows a comparison for the second case. In this example, $42$ constraints are active at the optimal solution. ALADIN just achieves a global linear convergence and is only three times faster than ADMM. 

\section{Conclusion}
This paper analyzed an optimization problem for coordinating multiple commercial buildings. The problem balances the voltage surge of the building network by using an energy cost defined demand response criterion. By introducing an auxiliary variable, the problem was reformulated into a standard distributed form with decoupled PWQ objectives and coupled affine equality constraints. For solving this non-smooth convex problem in an MPC scheme, we proposed a tailored ALADIN method, which can warmly start online and thus its convergence can be sped up. Our numerical results illustrated the effectiveness of the warm-start strategy and show that the ALADIN based MPC controller outperforms the ADMM based controller.

\section*{Acknowledgements}
{\small
JS, YJ and BH were supported by ShanghaiTech University under Grant-Nr. F-0203-14-012. CJ was supported from the Swiss National Science Foundation under the RISK project (Risk Aware Data Driven Demand Response, Grant-Nr. 200021 175627).}

\bibliography{ifacconf}

\end{document}